\documentclass[12pt, reqno]{amsart}
\usepackage{euscript, amssymb, amsmath, amsthm}
\usepackage{epsfig}
\usepackage{graphicx}

\usepackage{amssymb}
\usepackage{amsthm}

\usepackage{dsfont}
\usepackage[colorlinks=true,linkcolor=red]{hyperref}
\usepackage{lineno}
\usepackage{algorithm}
 \usepackage{amsmath} 
 \newtheorem{theorem}{Theorem}
 \newtheorem{remark}{Remark}
 \newtheorem{definition}{Definition}
 

\begin{document}
\title[PINNs for Timoshenko systems]{Physics-informed neural networks for Timoshenko system with Thermoelasticity}

\author[Chebbi et al.]{Sabrine Chebbi$^{1,*}$, Joseph Muthui Wacira$^{2}$, Makram Hamouda$^{3}$, Bubacarr Bah$^{4}$}

\address{$^{1}$African Institute for Mathematical Sciences, South Africa,
E-mail: sabrine@aims.ac.za\\ 
$^{2}$ African Institute for Mathematical Sciences, South Africa,  E-mail:joemuthui@aims.ac.za,\\
$^{3}$ Imam Abdulrahman Bin Faisal University, Dammam, Saudi Arabia, E-mail: mmhamouda@iau.edu.sa,\\ 
$^{4}$Medical Research Council Unit The Gambia at LSHTM, E-mail: bubacarr.bah1@lshtm.ac.uk
}

\maketitle

\begin{abstract}
The main focus of this paper is  to analyze the behavior of a numerical solution of the Timoshenko system coupled with Thermoelasticity and  incorporating second sound effects. In order to address this target, we employ the Physics-Informed Neural Networks  (PINNs) framework to derive an approximate solution for the system. Our investigation delves into the extent to which this approximate solution can accurately capture the asymptotic behavior of the discrete energy, contingent upon the stability number $\chi$. 
Interestingly, the PINNs overcome the major difficulties encountered while using the standard numerical methods.
\end{abstract}





\section{Introduction}\label{sec1}
Over the past decades, deep learning has revolutionized various fields by enabling computers to learn and extract complex patterns from vast datasets autonomously. Its utilization spans diverse domains, including image and speech recognition, natural language processing, autonomous vehicles, and medical diagnostics. Deep learning models, particularly neural networks with multiple layers, have shown their capacity to solve partial differential equations (PDEs) \cite{Rassi2}. To achieve this goal, the PDE solution is created via approximations using deep neural networks, and the partial derivatives are calculated through automated differentiation \cite{AUDif,12}. The formulation of a loss function is characterized by a mean-squared error norm (MSE), which includes the differential equations, along with the initial and boundary conditions. Minimizing this loss function across sampled points yields an approximate solution to the studied problem. This method is now widely recognized as PINNs, with initial investigations in \cite{23,22,21,B1, pinns}. This approach has been used to solve various PDEs, including elliptic, parabolic, hyperbolic, and inverse problems \cite{Foias,a, c}.

Despite significant advancements achieved in numerical techniques, including finite difference, finite element, finite volume, and spectral methods, persistent challenges are associated with the use of conventional numerical approaches. These challenges become particularly evident in scenarios involving multi-physics and multi-scale problems, where traditional numerical methods struggle to deliver dependable solutions (nonlinear terms, coupling, etc.).

The efficiency of PINNs and whether they can outperform those traditional numerical approaches have been addressed in \cite{Beat}, where both approaches are utilized to computationally solve a range of linear and nonlinear partial differential equations. These include Poisson equations in 1D, 2D, and 3D, Allen–Cahn equations in 1D, and semilinear Schrödinger equations in 1D and 2D. Subsequently, a comparison of computational expenses and the accuracy of approximations is conducted. 

Our primary goal in this article is to illustrate the applications of PINNs for complex coupled hyperbolic systems, both linearly and non-linearly damped, and analyze the asymptotic behavior of the energy, as time approaches infinity with respect to the coefficients of the system. To this end, we consider the Timoshenko system with thermoelasticity, incorporating second sound. This system is not only recognized as one of the fundamental models for depicting the lateral oscillation of a beam with fixed ends but also accounts for heat flux and temperature with respect to the Cattaneo law
$$\tau q_t + \beta q + \theta_x = 0,$$
where $q$ is the heat flux and $\beta$ is the coefficient of the thermal conductivity \cite{Santos}.  
Consequently,  we have the following system 
\begin{equation} 
\left\{
\begin{array}{ll}
\rho _{1}\varphi _{tt}-k(\varphi _{x}+\psi )_{x}=0,\  
 &  \ (x,t) \in  (0,L)\times (0,T), \\
\rho _{2}\psi _{tt}-b\psi _{xx}+k(\varphi _{x}+\psi )+\delta \theta
_{x}+\mu \psi_t=0,\ &  \ (x,t)\in \ (0,L)\times (0,T), \\
\rho _{3}\theta _{t}+q_{x}+\delta \psi _{xt}=0,\ &  \ (x,t)\in
(0,L)\times (0,T), \\
\tau q_{t}+\beta q+\theta _{x}=0,\ &    \ (x,t)\in \ (0,L)\times
(0,T).%
\end{array}%
\right.  \label{1}
\end{equation}%
where $t\in (0,T)$, with $T<T_{max}$, is the time variable and $0<T_{max} \le \infty$ is the maximal time of existence, $x\in (0,L)$ is the
space variable, the function $\varphi$ is the displacement vector, $\psi$
is the rotation angle of the filament, the function $\theta $ is the
temperature difference, $q=q(x,t)\in \mathbb{R}$ is the heat flux, and $\rho
_{1}$, $\rho _{2}$, $\rho _{3}$, $b$, $k$, $\delta $, $\beta $ and $L$ are positive
constants related to the physical properties of the beam. 

In the remainder of this article, we will assume without loss of generality that $L=1$.

Note that the main focus of this work is concerned with the obtaining of an accurate discrete energy to \eqref{1} in the context  a linear damping $\mu \psi_t$. Nevertheless, some additional numerical results are shown for other nonlinear dampings. Indeed, the discrete energy, as well as the theoretical exact one, is not depending on the damping term as we will see in \eqref{E} below. In view of \eqref{d} below, the non-increasing behavior of the energy is of course depending on the damping profile. 

We associate with \eqref{1} the following  mixed Dirichlet and Neumann boundary conditions
\begin{equation}
\varphi_x (0,t)=\varphi_x (1,t)=\psi (0,t)=\psi (1,t)=q(0,t)=q(1,t)=0,\hspace{%
0.58cm} \mbox{ }\forall \ t\in [0,T].  \label{cb}
\end{equation}%
Moreover, the initial conditions for the system \eqref{1}  are given by :
\begin{equation}
\left\{
\begin{array}{ll}
\varphi (x,0)=\varphi _{0}(x),\mbox{ }\varphi _{t}(x,0)=\varphi _{1}(x),\  
 &  \ x\in \ (0,1), \\
\psi (x,0)=\psi _{0}(x),\mbox{ }\psi _{t}(x,0)=\psi _{1}(x),\  
 &  \ x\in\  (0,1), \\
\theta (x,0)=\theta _{0}(x),\mbox{ }q(x,0)=q_{0}(x),\  
 &  \ x\in \ (0,1).%
\end{array}%
\right. \label{ci}
\end{equation}%
The aforementioned Timoshenko model has a wide range of applications in engineering, encompassing a diverse set of scenarios. These scenarios include tasks such as analyzing the behavior of suspension bridge cables, which often involve complex beam systems described by the Timoshenko beam equations \cite {h}. Additionally, Timoshenko models are relevant in contexts such as understanding the dynamic interactions between catenaries and pantographs in railway systems, where damped beam equations come into play \cite{g}. They also play a crucial role in simulating phenomena such as the impact of air turbulence on flight, which is essentialy captured by the Navier-Stokes equations \cite{Foias}, \cite{z}, and they apply to a multitude of other situations as well \cite{f}–\cite{S}.
In reality, the computation of analytical solutions of these equations is too complex, necessitating numerical approaches \cite{3}. Numerical techniques, such as finite-difference and finite-element methods, have been employed to approximate solutions for these PDEs. Although these methods have shown effectiveness, they encounter certain challenges. One such challenge is mesh generation, which becomes notably more complex when dealing with intricate geometries in higher dimensions. In \cite{2}, the authors employed a fourth-order finite difference scheme to calculate numerical solutions for the Timoshenko system coupled with thermoelasticity and second sound. This system, coupled with the Cattaneo law, resulted in four equations. The approach used by \cite{2} was adapted from the method utilized by \cite{Rap} 
yielding the decay rate for discrete solutions.
Note that  a few investigations has been conducted regarding the behavior and temporal evolution of the approximated energy. In this regard, in \cite{sab},  the decay rate of the discretized energy associated with the approximate solution of the 1D  Timoshenko system coupled by two equations, was derived using the Finite Element Method (FEM); see also \cite{MH} for an additional overview on the problem.
The existence and smoothness of solutions for the Timoshenko system incorporating distinct damping terms have been extensively explored in various works in the literature see \cite{albo, Fernd, Gass,  Riv}, along with the references therein.
We recall the well-posedness result to \eqref{1} in the next theorem for which the following functions' spaces are needed,
\begin{equation}
\begin{array}{c}
\displaystyle 
 L^2_{\star}(0,1):=\left\{ 
 v\in L^2(0,1)\  ; \ \displaystyle \int_{0}^1 v(s) ds=0 \right\}, \vspace{.2cm}\\ 
 \displaystyle 
 H^1_{\star}(0,1):=H^1(0,1) \cap L^2_{\star}(0,1),\vspace{.2cm}\\ \displaystyle 
 H^2_{\star}(0,1):=\left\{v\in H^2(0,1), \ v_x(0)=v_x(1)=0 \right\}.
\end{array}%
\end{equation}%
\begin{theorem}\textnormal{(\cite{2})} \label{thp}
 For all
initial data 
\begin{equation}
\begin{array}{ll}
    (\varphi_0, \varphi_1,\psi_0,\psi_1, \theta_0,q_0)& \in (H^2_{\star}(0,1) \cap H_{\star}^{1}(0,1) ) \times H^1_{\star}(0,1)\times (H^2(0,1)\cap H^1_0(0,1)) \nonumber\\
    &\times  H_{0}^{1}(0,1)\times H^1_{\star}(0,1)\times H^1_0(0,1), \nonumber
    \end{array}
   \end{equation}
    the system (\ref{1}), \eqref{cb} and \eqref{ci} has a unique
solution 
\begin{align}
\begin{array}{ll}
   (\varphi,\psi)& \in \mathcal{C}([0,T]; (H^2_{\star}(0,1) \cap H_{\star}^{1}(0,1) )\times (H^2(0,1)\cap H^1_0(0,1))  \cap \mathcal{C}^1([0,T];H^1_{\star}(0,1) \cap H_{0}^{1}(0,1)) \nonumber\\
    &\cap \mathcal{C}^2([0,T]; L^2_{\star} (0,1)\times L^2(0,1)), \nonumber
    \end{array}
  \end{align}
and 
$$(\theta,q) \in \mathcal{C}^0([0,T]; H^1_{\star}(0,1) \times H^1_0(0,1) \cap \mathcal{C}^1([0,T]; L^2_{\star}(0,1)\times L^2(0,1)).$$
\end{theorem}
\noindent Moreover, the energy associated with the solution denoted by
\begin{equation}\label{E}
 E(t)=E\left(\varphi,\psi,\theta,q\right)(t):= \frac{1}{2} \int_0^1\left(\rho_1\varphi_t^2+\rho_2\psi_t^2+b\psi_x^2+k(\varphi_x+\psi )^2+\rho_3 \theta^2+\tau q^2\right)dx.  
\end{equation}
 satisfies the following dissipation inequality 
\begin{equation} \label{d}
E^{\prime
}(t)=- \int^{1}_{0} \beta q^{2}dx-\int^{1}_{0} \mu \psi_{t}^2 dx \leq 0.
\end{equation} 
System \eqref{1} satisfies the strong stability result \eqref{stab} below, which has been shown in \cite{1} for a more general context of damping.  We  first consider the following conservative Timoshenko  system:
\begin{equation}\label{conservation}
\left\{
\begin{array}{ll}
\rho _{1}\varphi_{tt}-k(\varphi _{x}+\psi )_{x}=0,\quad & (x,t) \in  \
(0,1)\times (0,T), \\
\rho _{2}\psi_{tt}-b\psi _{xx}+k(\varphi _{x}+\psi )=0,\quad &   (x,t) \in  \ (0,1)\times (0,T).
\end{array}%
\right.
\end{equation}
Then, we assume the assumption  below  on the subset $\emptyset \neq\Omega\subset 
(0,1)$,
\begin{equation*}
(HS)\left\{
\begin{array}{ll}
\mbox{Let}\ (\varphi ,\psi)\
\mbox{ be a weak solution of \eqref{conservation}}  \\
\mbox{if}\ \psi _{t}\equiv 0\ \mbox{on}\ \Omega \ %
\mbox{then}\ (\varphi ,\psi )\equiv (0,0),%
\end{array}%
\right.
\end{equation*}%
and we formulate the stability result for the energy of %
\eqref{1} as follows.
\begin{theorem}\textnormal{(\cite{1})}
\label{TH2} Assume the hypothesis   $(HS)$ holds true. Then, for all 
\begin{align*}
    (\varphi_0, \varphi_1,\psi_0,\psi_1, \theta_0,q_0)& \in (H^2_{\star}(0,1) \cap H_{\star}^{1}(0,1) ) \times H^1_{\star}(0,1)\times (H^2(0,1)\cap H^1_0(0,1)) \nonumber\\
    &\times  H_{0}^{1}(0,1)\times H^1_{\star}(0,1)\times H^1_0(0,1), 
    \end{align*}
 the energy $E(t)$, defined by \eqref{E} and corresponding to the solution of %
\eqref{1}, satisfies
\begin{equation}  \label{stab}
 {\lim_{t\rightarrow \infty}E(t)=E_{\infty},}
\end{equation}
where $E_{\infty}=E(0,0,\theta_0,0)(0)$.
\end{theorem}
\begin{remark}
    \textnormal{In \cite{1}, the Timoshenko system coupled with thermoelasticity and second sound experiences indirect damping due to temperature coupling. Furthermore, the theorem presented above substantiates our proof that the system's energy converges to a non-zero equilibrium state $E_{\infty}$.}
\end{remark}
The subsequent sections of the article are structured as follows. In Section \ref{sec2}, we detail how PINNs is employed to approximate the solution of the Timoshenko system. In Section \ref{sec3}, we present the evolution in time of the PINNs approximate solution and the asymptotic behavior of the discrete energy with respect to the stability number 
\begin{equation}
    \label{xi-def}
\chi =\left( \tau -\frac{\rho _{1}}{k\rho _{3}}\right) \left( \frac{\rho _{2}}{b}-%
\frac{\rho _{1}b}{k}\right) - \frac{\tau \delta^{2}\rho _{1}}{{b k\rho _{3}}},
\end{equation}
for which it has been shown that an exponential decay rate  is obtained if and only if $\chi= 0$. Other types of damping are considered as well. Finally, in Section \ref{sec4}, we conclude by providing remarks and potential future directions.
\section{Physics Informed 
neural networks (PINNs) for Timoshenko system} \label{sec2}
 To solve the proposed Timoshenko system with Thermoelasticity \eqref{1}, we use the following neural network as an approximate solution 
\begin{eqnarray}
\begin{array}{ll}
    \varphi(x,t)\simeq \Tilde{\varphi}(x,t)=\varphi_{\Xi}(x,t,\Xi),  \
      \psi(x,t)\simeq \Tilde{\psi}(x,t)=\psi_{\Xi}(x,t;\Xi),\nonumber \vspace{.2cm}\\
   \theta(x,t) \simeq \Tilde{\theta}(x,t)=\theta_{\Xi}(x,t;\Xi), \ \ 
     q(x,t) \simeq \Tilde{q}(x,t)=q_{\Xi}(x,t;\Xi).
     \end{array}
\end{eqnarray}
where $\Xi=\{W,b\}$ is the networks' weights and biases tuple.

We define the residuals associated with the Timoshenko system \eqref{1} as follows:
\begin{equation}
 \label{RPDE}
\begin{array}{ll}
  &\mathcal{R}_{PDE_1}([\Tilde{\varphi},\Tilde{\psi}])(x,t)=\rho _{1}\Tilde{\varphi}_{tt}-k(\Tilde{\varphi} _{ x}+\Tilde{\psi})_{x}, \vspace{.2cm}\\
 &\mathcal{R}_{PDE_2}([\Tilde{\psi},\Tilde{\varphi},\Tilde{\theta}])(x,t)= \rho _{2}\Tilde{\psi} _{tt}-b\Tilde{\psi} _{xx}+k(\Tilde{\varphi} _{x}+\Tilde{\psi} )+\delta \Tilde{\theta
}_{x}+\mu \Tilde{\psi}_t, \hspace{0.2cm}  \vspace{.2cm}\\
&\mathcal{R}_{PDE_3}([\Tilde{\theta},\Tilde{q},\Tilde{\psi}])(x,t)=\rho _{3}\Tilde{\theta} _{t}+\Tilde{q}_{x}+\delta \Tilde{\psi} _{xt},  \vspace{.2cm}\\ 
 &\mathcal{R}_{PDE_4}([\Tilde{q},\Tilde{\theta}])(x,t)=\tau \Tilde{q}_{t}+\beta 
 \Tilde{q}+\Tilde{\theta} _{x}.
 \end{array}
  \end{equation}
  \\
Additionally, we define the residuals associated with the  boundary conditions and initial conditions, \eqref{cb} and \eqref{ci}, as shown in \eqref{Rbc} and \eqref{I_cond-pinn} below, respectively. These residuals will be used to train the neural network.
\begin{itemize}
    \item  The mixed  Dirichlet and Neumann Boundary conditions for all $t\in[0,T]$:
    
\begin{equation}\label{Rbc}
\begin{array}{ll}
\mathcal{R}_0^b([\Tilde{\varphi},\Tilde{\psi},\Tilde{q}])(0,t)=\Tilde{\varphi}_x(0,t)-g_1(t)+\Tilde{\psi}(0,t)-g_2(t)+ \Tilde{q}(1,t)-g_3(t), \vspace{.2cm}\\
\mathcal{R}_1^b([\Tilde{\varphi},\Tilde{\psi},\Tilde{q}])(1,t)= \Tilde{\varphi}_x(1,t)-\hat{g}_1(t) +\Tilde{\psi}(1,t)-\hat{g}_2(t)+ \Tilde{q}(1,t)-\hat{g}_3(t),
\end{array}
\end{equation} 
  \item  The Initial conditions for all $x\in (0,1)$:
\begin{equation}
\begin{array}{rcl}
\mathcal{R}_{1}^I([\Tilde{\varphi},\Tilde{\psi},\Tilde{\theta},\Tilde{q}])(x,0) &=& \Tilde{\varphi}(x,0)-\varphi_0(x) + \Tilde{\psi}(x,0)-\psi_0(x) \vspace{.2cm}\\ &&+\Tilde{\theta}(0,x)-\Tilde{\theta}_0(x)+\Tilde{q}(0,x)-\Tilde{q}(x), \vspace{.2cm}\\
\mathcal{R}_{2}^I([\Tilde{\varphi},\Tilde{\psi}])(x,0) &=& \Tilde{\varphi}_t(x,0)-\varphi_1(x)+\Tilde{\psi}_{t}(x,0)-\psi_1(x) .
   \label{I_cond-pinn}
\end{array}
\end{equation}
\end{itemize}
Note that the exact solution $(\varphi,\psi,\theta,q)$ satisfies  $$ \mathcal{R}_{PDE_1}([\varphi,\psi])= \mathcal{R}_{PDE_2}([\varphi,\psi,\theta])= \mathcal{R}_{PDE_3}([\theta,q,\psi])= \mathcal{R}_{PDE_4}([\theta,q])=0,$$
and
$$\mathcal{R}_{2}^I([\varphi,\psi])=\mathcal{R}_{1}^I([\varphi,\psi,\theta,q])=\mathcal{R}_1^b([\varphi,\psi,q])=\mathcal{R}_1^0([\varphi,\psi,q])=0.$$
Let $S$ be a set of  collocation points $\{X=(x_i,t_j) \}_{i,j=1}^{N}, N\geq 3000$. We define \textit{the training loss} function, called also \textit{loss function} as
\begin{equation} \label{LOSS}
\mathcal{E}(\Xi,S)^2:= \sum_{k=1}^8 MSE_{k},
\end{equation} where the Mean Squared Error  $MSE_{k}$ for $k\in\{1,\dots,8\}$ is given by
\begin{eqnarray}
    MSE_{1}=\frac{1}{N}\sum_{i,j=1}^{N}\left|\mathcal{R}_{PDE_1}([\Tilde{\varphi},\Tilde{\psi}])(x_i,t_j) \right|^2, \hspace{0.3cm} x_i\in(0,1), \ t_j\in (0,T), \nonumber\\
    MSE_{2}=\frac{1}{N}\sum_{i,j=1}^{N} \left|\mathcal{R}_{PDE_2}([\Tilde{\psi},\Tilde{\varphi},\Tilde{\theta}])(x_i,t_j)\right|^2, \hspace{0.3cm} x_i\in (0,1), \ t_j\in(0,T),\nonumber \\
     MSE_{3}=\frac{1}{N}\sum_{i,j=1}^{N}|\mathcal{R}_{PDE_3}([\Tilde{\psi},\Tilde{\varphi},\Tilde{\theta}])(x_i,t_j)|^2, \hspace{0.3cm} x_i\in (0,1), \ t_j\in(0,T),\nonumber \\
      MSE_{4}=\frac{1}{N}\sum_{i,j=1}^{N}\left|\mathcal{R}_{PDE_4}([\Tilde{\psi},\Tilde{\varphi},\Tilde{\theta}])(x_i,t_j)\right|^2, \hspace{0.3cm} x_i\in (0,1), \ t_j\in(0,T),\nonumber \\
     MSE_{5}=\frac{1}{N} \sum_{j=1}^{N}\left|  \mathcal{R}_0^b[\Tilde{\psi},\Tilde{\varphi},\Tilde{q}](0,t_j)\right|^2, \ \ 
      MSE_{6}=\frac{1}{N} \sum_{j=1}^{N}\left| \mathcal{R}_1^b[\Tilde{\psi},\Tilde{\varphi},\Tilde{q}](1,t_j)\right|^2,  \ t_j\in(0,T),\nonumber \\
        MSE_{7}=\frac{1}{N} \sum_{i=1}^{N} \left|\mathcal{R}_{1}^I([\Tilde{\varphi},\Tilde{\psi},\Tilde{\theta},\Tilde{q}])(x_i,0)\right|^2, \ \
        MSE_{8}=\frac{1}{N} \sum_{i=1}^{N}|\mathcal{R}_{2}^I[\Tilde{\varphi},\Tilde{\psi}](x_i,0)|^2, \hspace{0.3cm} x_i\in (0,1).\nonumber 
\end{eqnarray}
The terms of \eqref{LOSS} may be scaled by different weights (penalty coefficients). For simplicity, we set
all these weights to one in our analysis. Hence, PINNs aim to minimize the squared training error $\mathcal{E}(\Xi, S)^2$ with respect to the network parameters $\Xi$. Once the optimization process converges, the trained network $(\Tilde{\varphi}, \Tilde{\psi}, \Tilde{\theta}, \Tilde{q})$ provides an approximation for the solution $(\varphi, \psi, \theta, q)$ of the system \eqref{1}.
\section{Numerical experiments} \label{sec3}
In this section, we provide a comprehensive explanation of the code and the architecture employed in the PINN used to approximate the solution of the Timoshenko beam system.
\subsection{ Discussion about the neural network architecture}
\begin{figure}[!ht]
\includegraphics[width=14cm]{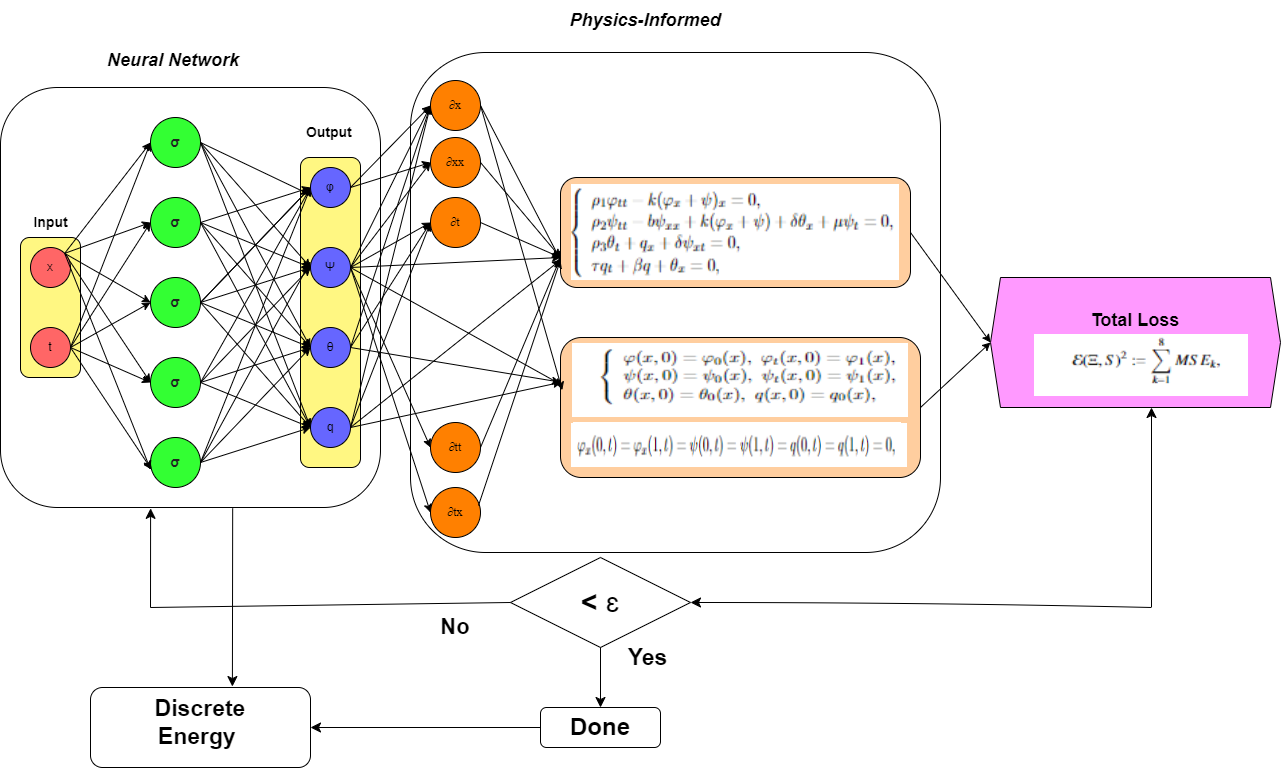}
\centering
\caption{Schematic diagram of the PINNs framework for Timoshenko systems.} 
\label{fig:kkg}
\end{figure}
We constructed a feedforward neural network using PyTorch for our study. The network features two input neurons, corresponding to time ($t$) and space ($x$), followed by five hidden layers with 100 neurons each. The output layer comprises four neurons that represent our system's dependent variables. We utilize the hyperbolic tangent (tanh) activation function for the hidden layers and a linear activation function for the output layer. The primary goal is to minimize the Mean Squared Error (MSE) defined by \eqref{LOSS}, and we achieve this using the Adam optimizer with a learning rate of 0.0005. Our training dataset consisted of 3000 randomly selected collocation points within the bounds of time and space. The network is trained for 10000 epochs on Google Colab with the free T4 GPU.
\subsection{PINNs for Timoshenko system with an exact solution}
There exist several methods that can be used to validate PINNs solutions, especially when analytical solutions are unavailable. One approach involves comparing PINNs results to those generated using numerical techniques such as finite difference, finite element, finite volume, or spectral methods. Another method entails comparing the PINNs output to experimental data. Here, the goal is to match the predicted solutions from PINNs with empirically measured values across both space and time, thereby establishing a connection between the model's predictions and real-world observations. In this section, we perform
some experiments to ensure the convergence of the PINNs solution to the exact one  in the sense that the relative error $\mathcal{R}$ is small enough. To this end, the following problem is considered.\\
Set 
$\rho_1=\rho_2=\rho_3=\tau=\delta=k=\beta=b=1, \  L=1$ and $\mu=1$ in \eqref{1}, this yields
\begin{equation}  \label{exact}
\left\{
\begin{array}{ll}
\varphi _{tt}-(\varphi _{x}+\psi )_{x}=f_1,\quad & \mbox{in}\ (0,1)\times (0,T), \\
\psi _{tt}-\psi _{xx}+(\varphi _{x}+\psi )+ \theta
_{x}+\psi_t=f_2,\quad & \mbox{in} \ (0,1)\times (0,T), \\
\theta _{t}+q_{x}+ \psi _{xt}=f_3,\quad & \mbox{in}\
(0,1)\times (0,T), \\
 q_{t}+ q+\theta _{x}=f_4,\quad & \mbox{in} \ (0,1)\times
(0,T).%
\end{array}%
\right. 
\end{equation}%
The Cauchy problem corresponding to \eqref{exact} is associated with explicit initial data, namely
\begin{equation} \label{exp-inti}
\left\{
\begin{array}{l}
    \varphi_{0}(x)=4 x(1-x)= \varphi_1(x), \\  \psi_{0}(x)= 4 x(1-x) = \psi_1(x), \\
      \theta_{0}(x)=4 x(1-x), \\  q_{0}(x)= 4 x(1-x),
\end{array}
\right. 
\end{equation}
together with the following boundary conditions:
\begin{align}
  \varphi (0,t)=\varphi(1,t)=\psi (0,t)=\psi (1,t)=q(0,t)=q(1,t)=0,\hspace{%
0.58cm}\forall \mbox{ }t\in [0,T].  \label{Mcb}
\end{align}
Note that, in \eqref{exact}, we denoted by $f_1,\ldots, f_4$ the source terms  given as follows:
\begin{equation} \label{St}
\left\{
\begin{array}{l}
    f_1(x,t)=4e^t(-x^2+3x+1), \\ f_2(x,t)=4e^t(-3x^2-x+4),\\
    f_3(x,t)= 4e^t(-x^2-3x+2), \\ 
    f_4(x,t)=4e^t(-2x^2+1).
\end{array}
\right. 
\end{equation}
The exact solution to the problem \eqref{exact}, taking into account the   explicit expressions of the sources terms \eqref{St}, is given by 
\begin{equation} \label{sext}
\left\{
\begin{array}{l}
    \varphi_{ex}(x,t)=4e^t x(1-x), \\ \psi_{ex}(x,t)= 4e^t x(1-x),\\
      \theta_{ex}(x,t)=4e^t x(1-x), \\ q_{ex}(x,t)= 4e^t x(1-x),
\end{array}
\right. 
\end{equation}
which is consistent with the boundary conditions \eqref{Mcb}.

\begin{figure}[!ht]
\includegraphics[width=6.7cm]{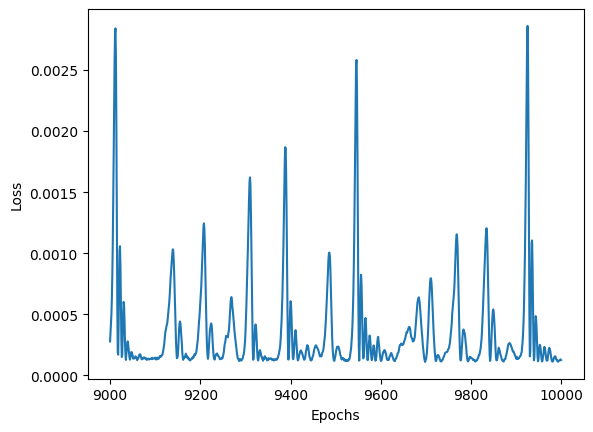}
\includegraphics[width=5.5cm]{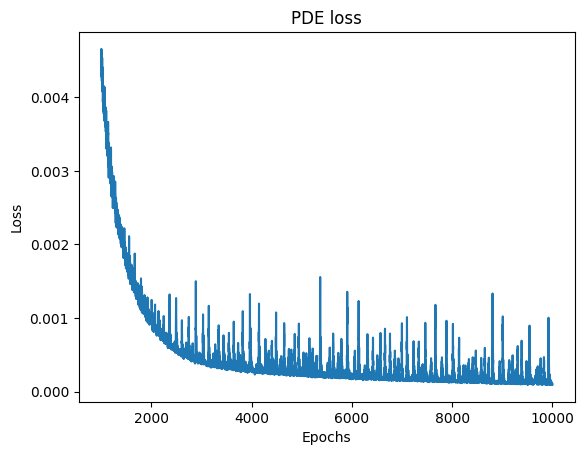}
\centering
\caption{The value of the total loss and PDE loss as a function of the number of epochs}
\label{fig:Tloss}
\includegraphics[width=6cm]{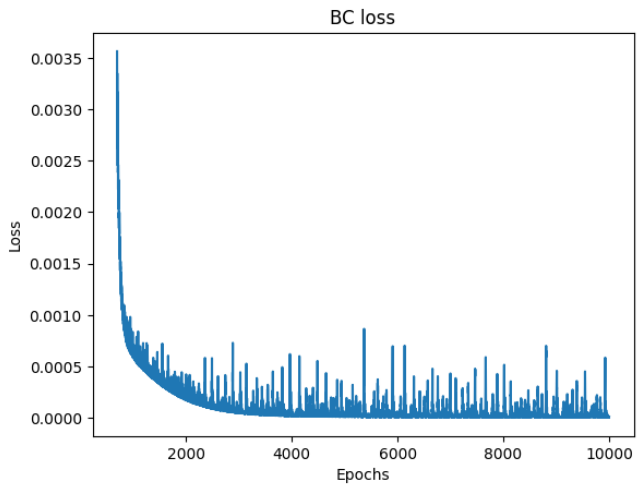}
\includegraphics[width=6cm]{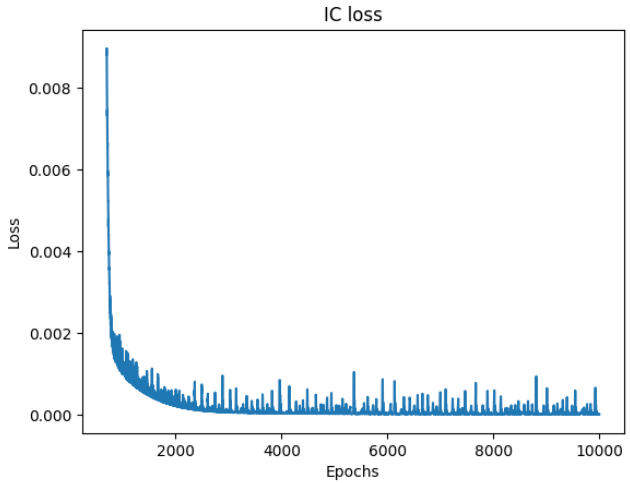}
\centering
\caption{The  boundary and the initial conditions losses.}
\label{fig:00}
\includegraphics[width=6.5cm]{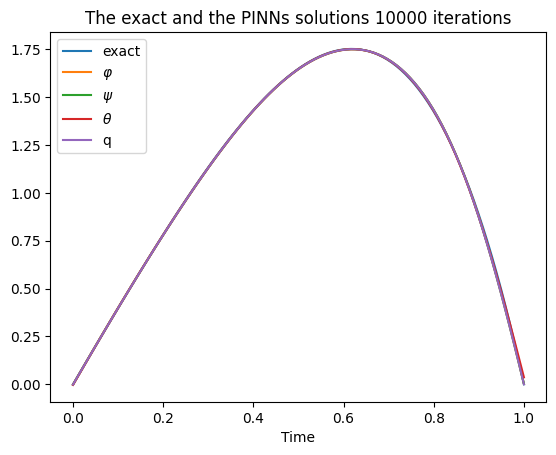}
\includegraphics[width=5.5cm]{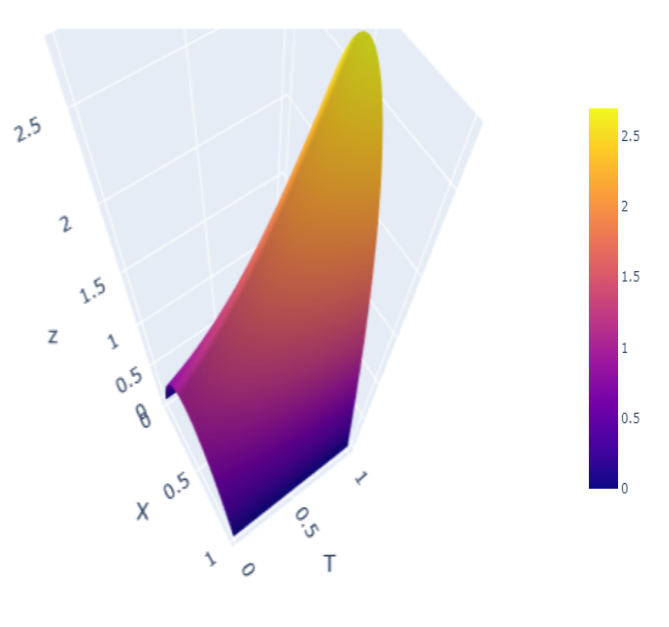}
\centering
\caption{PINNs predicted solutions $(\varphi,\psi,\theta,q)$  superimposed on the exact solution}
\label{fig:kk}
\end{figure}

We start by implementing the exact solution \eqref{sext} for the linear damped Timoshenko system \eqref{exact}. This choice of solution ensures that the initial conditions and the Dirichlet boundary conditions match with those outlined in \eqref{exp-inti} and \eqref{Mcb}. The behavior of the total loss function and the PDE loss are shown in Figure \ref{fig:Tloss}. Figure \ref{fig:00} exhibits a comprehensive view of the individual losses that make up the total loss {\it i.e.} losses associated with  the boundary conditions and initial conditions. 
As shown in Figure \ref{fig:00}, the network undergoes training for $10^4$ epochs. The aforementioned figure clearly demonstrates that the optimizer converges to the exact solution, resulting in an upper training loss of $\approx 9 \times 10^{-3}$ for the initial conditions and around $3.7\times 10^{-2}$ for the boundary conditions loss. The graph reveals that the optimizer approaches a local minimum of nearly zero. Therefore, increasing the number of epochs with the same neural network configuration will likely lead to further convergence.
\begin{figure}[!ht]
\includegraphics[width=6.5cm]{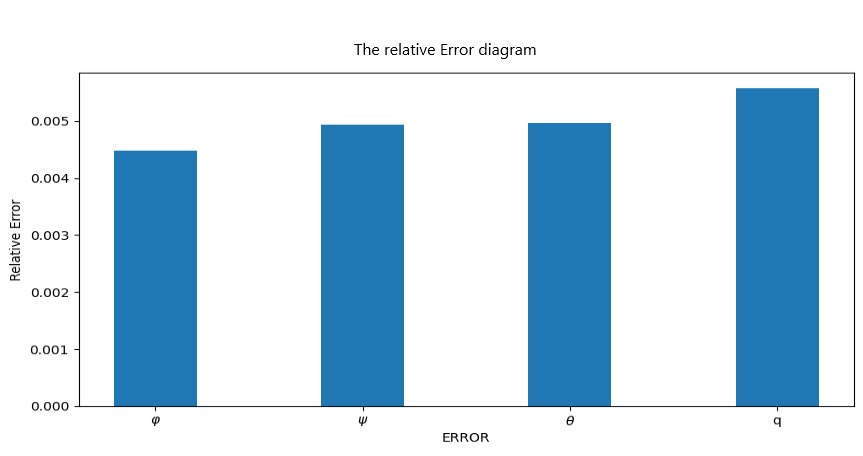}
\includegraphics[width=6.5cm]{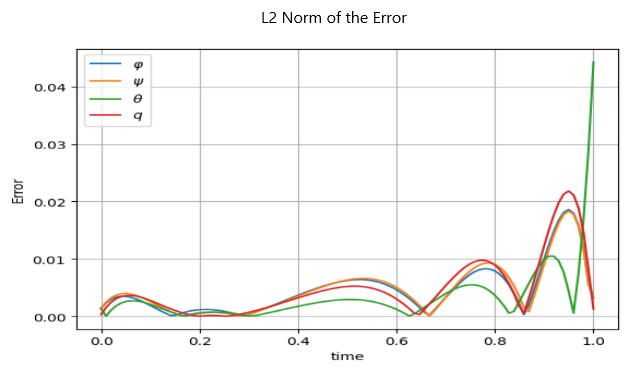}
\centering
\caption{\textbf{Left:}Thermoelastic Timoshenko beam  $\mathbb{L}^2$ norm error in predictions; \textbf{Right:} Relative error }
\label{fig:kk1}
\end{figure}
  Figure \ref{fig:kk} shows the exact solution alongside the PINNs solution. We observe how closely the two solutions align in behavior. To confirm the fact that the PINNs solution exhibits a similar behavior to the exact solution, we plot the $\mathbb{L}^2-$error shown in Figure \ref{fig:kk1}. For  the error estimation,  the relative error $\mathcal{R}$ is also used as in \cite{INVTiMO} but without  percentage. Here, $\tilde{u}$
is the prediction and $u$ is the
contracted exact solution.
 We define the Relative error $\mathcal{R}$ by 
  \begin{equation}
      \mathcal{R}= \frac{\|\tilde{u}-u\|_2}{\|u\|_2},
  \end{equation}
  The Relative error corresponding to each component of the solution $(\varphi, \psi, \theta, q)$ is presented in Figure \ref{fig:kk1}. 
\begin{remark}
   \textnormal{As evident in Figure \ref{fig:kk1}, the $\mathbb{L}^2-$error illustrates a deviation in the error associated with the temperature variable $\theta$, reaching a relatively higher value of approximately $5 \times 10^{-2}$. This discrepancy arises because Physics-Informed Neural Networks (PINNs) necessitate explicit information, particularly for the boundary condition related to $\theta$, see also Figure \ref{fig:00} for the BC loss. However, theoretically, the boundary conditions for $\theta $ are implicitly derived from the other unknowns $(\varphi,\psi,q)$ using the equations in \eqref{1} and the boundary conditions \eqref{cb}. Introducing this additional explicit information would theoretically over determine the system described by \eqref{1}--\eqref{ci}. We believe that this observation will improve the present results, but we postpone this discussion to a forthcoming work. }
\end{remark}
\subsection{PINNs on stability of Timoshenko system}\mbox{}\\
In the context of the Timoshenko system with a second sound without any damping term, a new number called $\chi$ has been introduced to quantify the exponential decay \cite{Stabn}. The system is considered to be exponentially stable if and only if $\chi$ is equal to zero. However, if $\chi$ is non-zero, the system lacks exponential stability, and the semigroup has a polynomial decay. 

\begin{itemize}
    \item \textbf{Case 1.} We consider the system \eqref{1} with $\mu=0$ and  we assume in this case that the stability number $\chi$ is non equal to zero. 
    Solving system \eqref{1} while considering the Dirichlet Boundary conditions for the case of simplicity, $$\varphi|_{x=0,x=1}=\psi|_{x=0,x=1}=q|_{x=0,x=1}=0,$$ using the PINNs framework. We consider the following parameters: the number of collocation points is $N_x=N_t=3000$, the final time is $T=40$, the number of iterations is $8000$, the damping coefficient $\mu=0$ (indicating no additional damping term), and $\chi \neq 0$ we can take different values of $\rho_1,\rho_2, \rho_3, k$ and $\tau $ as far as $\chi$ is always non zero.
    \item  \textbf{Case 2.} In this case, 
     we consider the system \eqref{1} for $\mu=0$ and  we choose the data so  that the stability number $\chi$ is equal to zero.  
We take $T=30$, iterations $=8000$, and we keep the same number of collocation points $N_x=N_t=3000$. More precisely, we set
$$\rho_1=\rho_2= k=2, \ 
\rho_3=b=\beta=1,\ 
\delta=\sqrt{\frac{2}{3}}, \ \tau=3. $$
\end{itemize}
\begin{figure}[!ht]
\includegraphics[width=6.5cm]{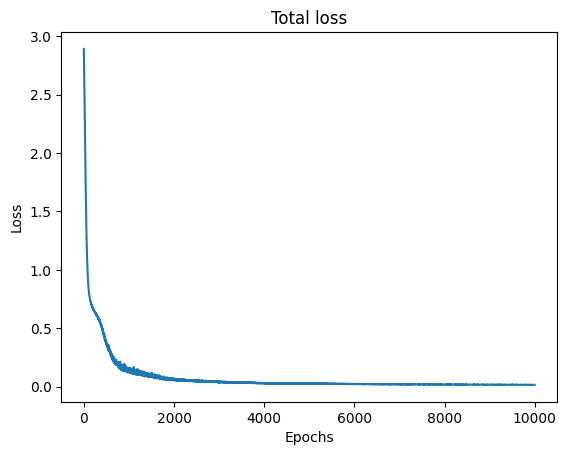}
\includegraphics[width=6.5cm]{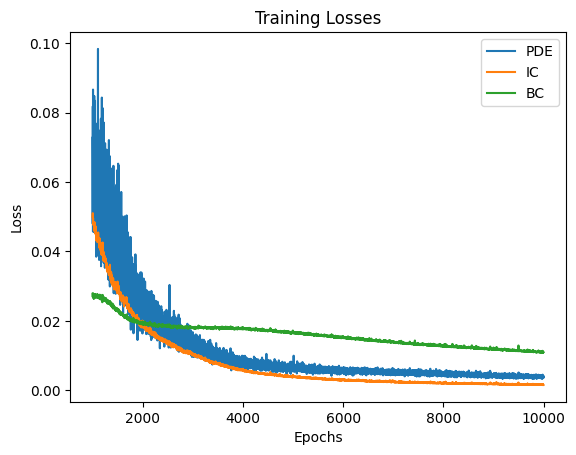}
\centering
\caption{\textbf{Left:}  the Total loss;  \textbf{right:} the PDE, the boundary and the initial conditions' losses for case 1}
\label{fig:3m}
\end{figure}
\begin{figure}[!ht]
\includegraphics[width=6.5cm]{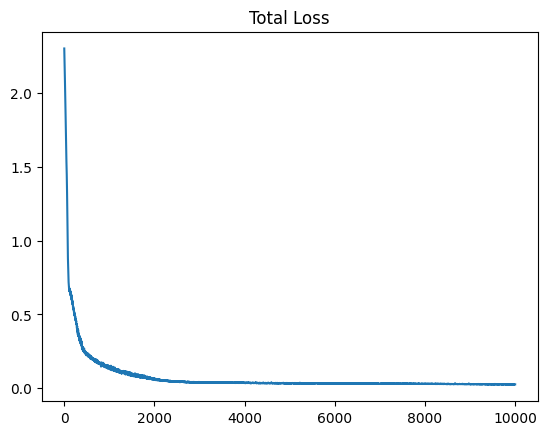}
\includegraphics[width=6.5cm]{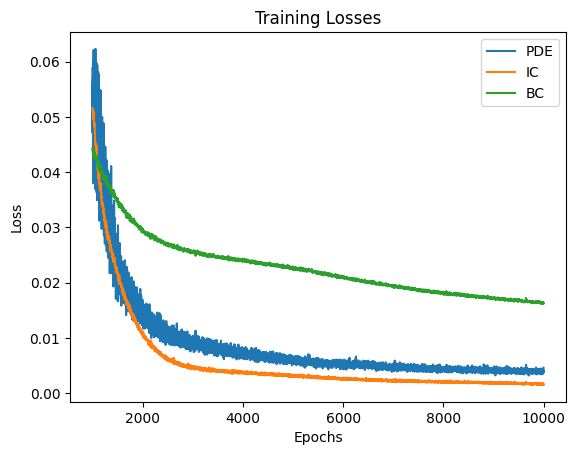}
\centering
\caption{\textbf{Left:} The Total loss; \textbf{Right:} The PDE boundary and the initial conditions' losses for case 2}
\label{fig:31m}
\end{figure}
Figures \ref{fig:3m} and  \ref{fig:31m} illustrate the behavior of the total loss functions corresponding to the different characteristics of the problem. The figures clearly demonstrate that the optimizer converges to the solution of the Timoshenko system \eqref{1} in both cases (case 1 and case 2), resulting in a training loss that reaches the value $\approx 5\times 10^{-1}$ for the total loss, a loss of order $\approx 1\times 10^{-1}$ for the PDE, a loss $\approx 5\times 10^{-1}$ for the initial condition  and  for the boundary condition the loss is of order $\approx 5\times 10^{-2}$. For instance,  for epochs in $[0,2000]$ for case 2, see right of Figure \ref{fig:31m} and the loss is about $\approx 2\times 10^{-2}$ for case 1 as shown in the right of Figure \ref{fig:3m}. The aforementioned values are in fact upper bounds for the losses when the number of epochs is a bit low. Nevertheless, the graph reveals that, upon increasing the number of epochs with the same neural network configuration, we  clearly obtain  a more accurate convergence.
\begin{figure}[!ht]
\includegraphics[width=6.5cm]{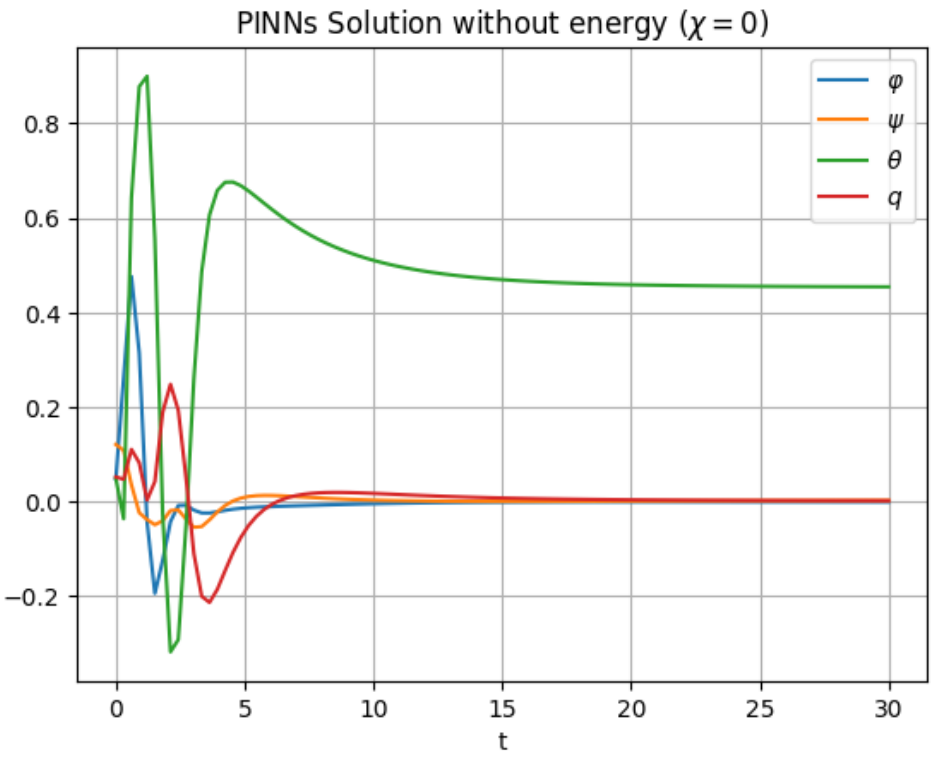}
\includegraphics[width=6.5cm]{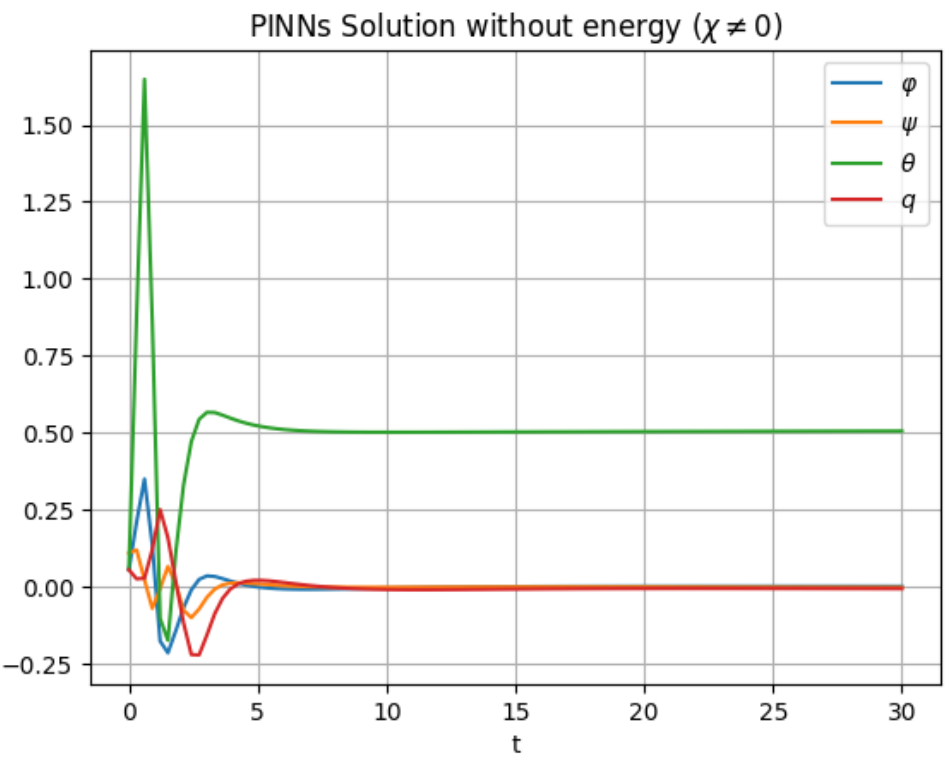}
\centering
\caption{The final time ($T=30$) profile  of the Solution $(\varphi,\psi,\theta,q)$ with $\mu=0$, zero-stability number $\chi = 0$ (left) and non-zero stability number $\chi \neq 0$ (right).  }
\label{fig:123}
\end{figure}
In Figure \ref{fig:123}, we present the solution of the Timoshenko system without adding any damping term, namely $\mu =0$.  The solution exhibits a clear decreasing over time, with the presence of some oscillations before reaching a stable state; this occurs starting from time $t\geq 5$ for  $(\varphi,\psi,q)$ which approach $(0,0,0)$. 
 However, for the temperature $\theta$, the limiting  value is less than $0.5$ for the case when $\chi=0$, while  the stable value of $\theta $ in the case 2 ($\chi\neq 0$) is $0.5$. The main interpretation of this figure is that the strong stability of the system is clearly justified, and it is independent of  the stability number $\chi$, and consequently   the physical parameters of the system. This confirms  the theoretical results in Theorem \ref{TH2} especially the  different stability regimes  for the cases $\chi=0$ and $\chi\neq 0$. Nevertheless, the strong stability of the energy remains consistent. 

\subsection{Energy decay rate analysis}
Recalling the expression of the energy associated with the solution to \eqref{1}, which is given by \eqref{E}, it is worth to mention that it would be not practical, in a general context,  to compute the integral in  \eqref{E}. Therefore,  in a natural way, we need to approximate the integral in the energy using the  Numerical Quadrature Rules. The basic idea of the quadrature rule is to replace the definite integral with a sum of the function values evaluated at certain points (called quadrature points) multiplied by  numbers (called
quadrature weights).
\begin{remark}
\textnormal{At this stage, we utilize the set of collocation points $(x_i, t_j)$ as inputs, along with the PINNs output, which represents the predicted Timoshenko solution denoted as $(\phi, \psi, \theta, q)$. We refer to Figure \ref{fig:123} for a more detailed illustration of its behavior over time. To approximate the derivatives with respect to both $x$ and $t$, we employ a Finite Difference  scheme.
We evaluate the energy associated with the PINNs solution after only $1000$ iterations and this is enough to analyze the asymptotic behavior of the discrete energy.}
\end{remark}
\begin{algorithm}[H]
\caption{Disecte energy for PINNs solution of Timoshenko system \eqref{1}}\label{alg:cap}
\begin{itemize}
    \item[(i)] Compute the PINNs solution 
\item[(ii)] Use Finite Difference method in time ($t$) and space ($x$) to compute the derivatives of the PINNs solution obtained in (i).
\item[(iii)] Evaluate the discrete energy inherited from the PINNs predicted solution.
\end{itemize}
\end{algorithm}
\begin{definition}(\cite{4})
  Let $\Lambda  \subset \mathbb{R}$ be the domain of the function $f\in L^1(\Lambda)$, then the quadratic rule formulae provides an approximation of $\int_{\Lambda} f(s) ds$  as 
\begin{equation}
   \int_{\Lambda} f(s) ds  \approx \frac{1}{M} \sum_{n=1}^M w_n f(s_n),
\end{equation}
where $s_n\in \Lambda \  (1\leq n\leq M)$   are the quadrature points and $w_n$ are the suitable
quadrature weights. The approximation accuracy is influenced by the type of quadrature rule, the number
of quadrature points (M), and the regularity of the  function $f$. 
\end{definition}
To approximate the integral in \eqref{E}, we set the time quadrature points in the dataset $(t_i)^{i=N_t}_{i=0}$, and let $(x_j)^{j=N_x}_{j=0}$ be the set of collocation points in $(0,1)$. Therefore, the approximated energy is given by 
\begin{eqnarray}\label{Eee}
    \mathcal{E}_i&=&\frac{1}{2 N_x} \sum_{j=0}^{N_x-1} \Big[
     \rho_1 \left|\frac{\tilde{\varphi}(x_{j+1},t_{i+1})-\Tilde{\varphi}(x_{j+1},t_{i})}{t_{i+1}-t_i}\right|^2 +\rho_2 \left|\frac{\tilde{\psi}(x_{j+1},t_{i+1})-\tilde{\psi}(x_{j+1},t_i)}{t_{i+1}-t_i}\right|^2 \nonumber\\&&+b\left|\frac{\tilde{\psi}(x_{j+1},t_i)-\tilde{\psi}(x_j,t_i)}{x_{j+1}-x_j}\right|^2
     +k\left|\frac{\Tilde{\varphi}(x_{j+1},t_i)-\tilde{\varphi}(x_j,t_i)}{x_{j+1}-x_j}+\tilde{\psi}(x_{j+1},t_i)\right|^2  \\&&+\rho_3 |\tilde{\theta}(x_{j+1},t_i)|^2+\tau |\tilde{q}(x_{j+1},t_i)|^2\Big],\nonumber
\end{eqnarray}
where $N_t$ and  $N_x>0$ are the number of collocation points in time and space, respectively.
\subsubsection{Non conservation of the energy,  $\mu=0 $.}
\begin{figure}[!ht]
\centering
\includegraphics[width=6.5cm]{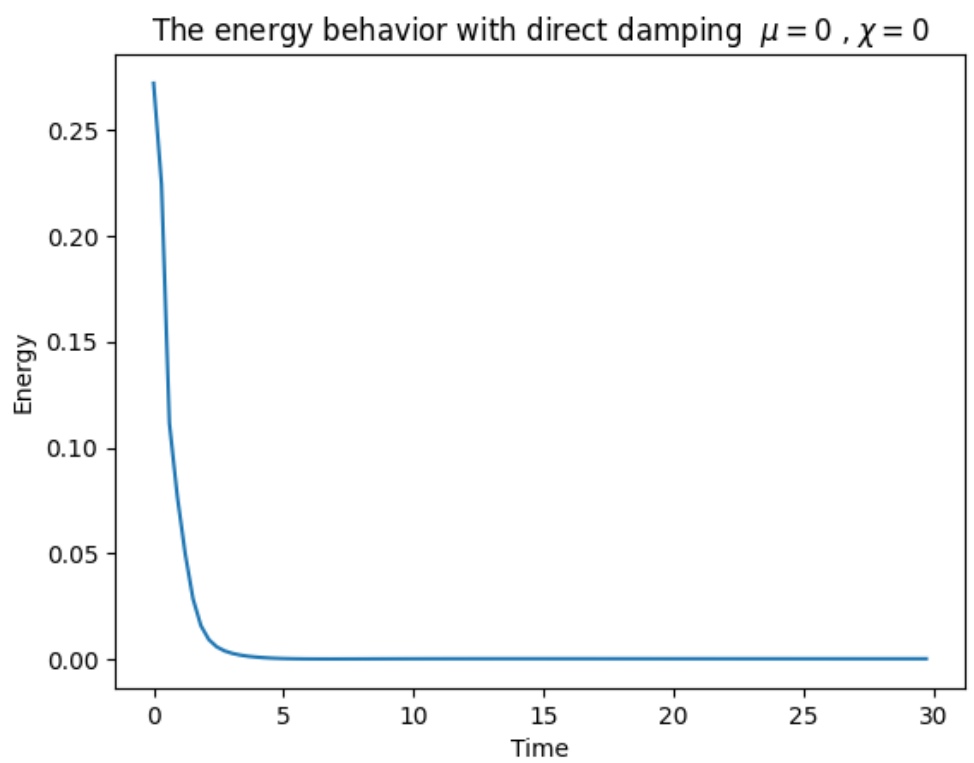}
\includegraphics[width=6.5cm]{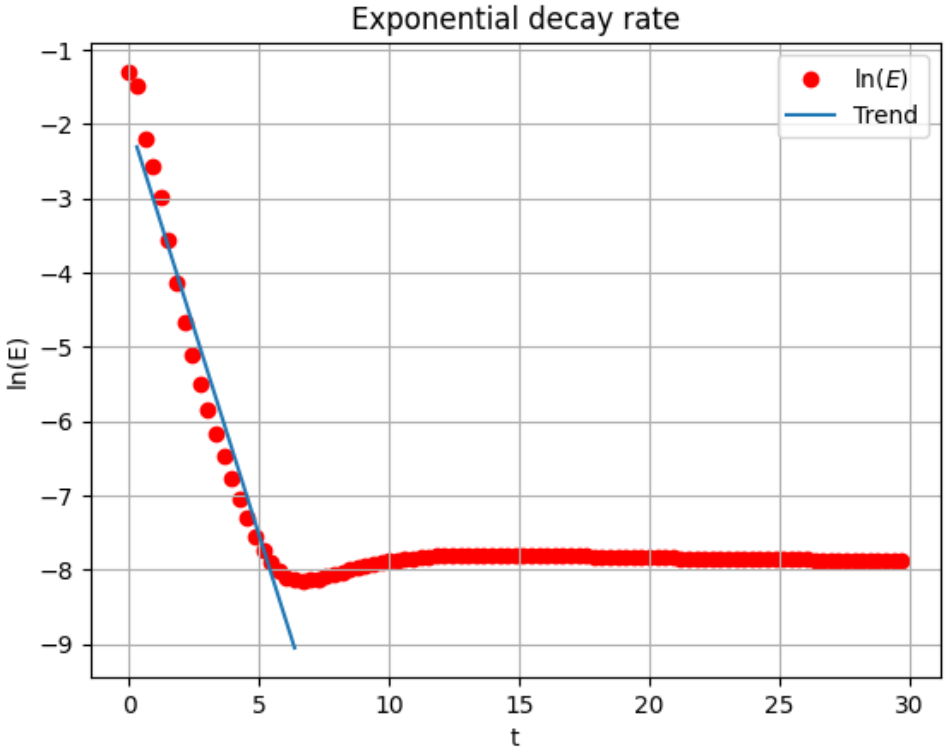}
\caption{The exponential behavior of the energy, $\mu=0$ and $\chi=0$.}
\label{fig:expo}
\end{figure}

Taking $\rho_1=\rho_2= k=2, \ 
\rho_3=b=\beta=1,\ 
\delta=\sqrt{\frac{2}{3}}, \ \tau=3 \text{\quad and\quad} \mu=0$, we plotted the energy as a function of time and observed that it has a decreasing curve, which means that $\frac{dE(t)}{dt}<0$. This proves the non-conservation of the energy associated with the solution of system even without adding any damping term. 
 In  Figure \ref{fig:expo},  we took $\theta_0(x)= 4x(1-x)$,   then as it is shown in the left-hand side of Figure \ref{fig:expo}, we observe that 
 $$\lim_{t \to \infty}E(t)= E_{\infty}= 0.$$ 
\begin{figure}[!ht]
\includegraphics[width=5.5cm]{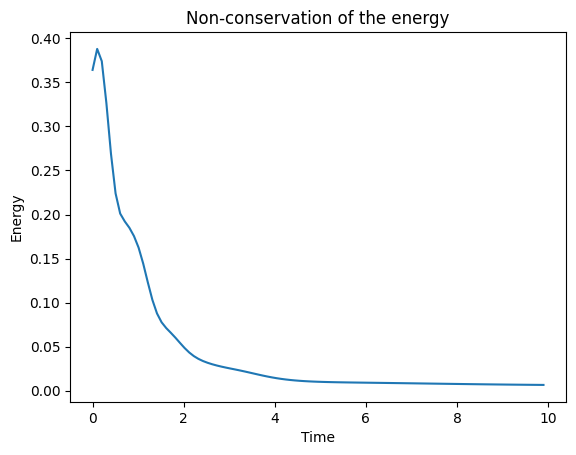}
\includegraphics[width=6.5cm]{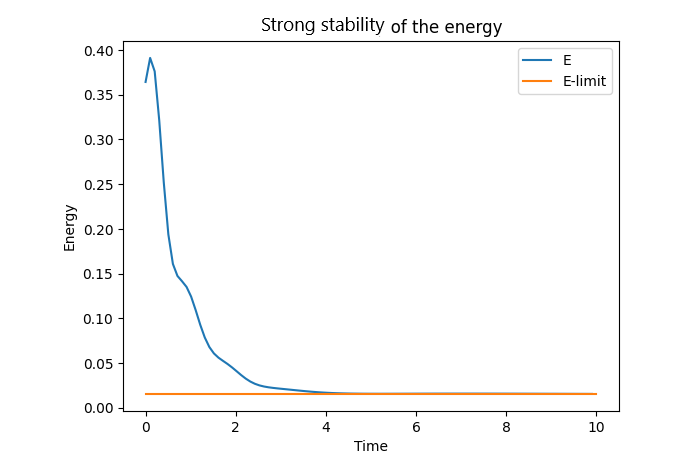}
\centering
\caption{The strong stability of the energy, $\mu=0$ and $\chi=0$, $\theta_0(0)\neq 0.$ (internal damping).}
\label{fig:d77}
\end{figure}
In Figure \ref{fig:d77}, we have illustrated the energy behavior for  $T=10$ while incorporating a non-zero initial value for $\theta_0(0)$ starting from $T=4$, the energy reaches a stable state with a value $\approx 2.10^{-2}$. Consequently, we noticed that the energy evolves from an initial level nearly equal to $0.4$, surpassing the initial value depicted in the left side of Figure \ref{fig:expo} at approximately $0.3$. Notably,  we clearly see the horizontal asymptotic line giving rise to the value of $E_{\infty} \approx 2.10^{-2}$. In agreement with the strong stability statement  in Theorem \ref{TH2}, we have that $\lim_{t \to \infty}E(t)= E_{\infty}$ non zero value.  
Figure \ref{fig:expo} has shown a consistently decreasing trend of the energy plotted against time. In the plot on the right side of  Figure \ref{fig:expo}, we plotted $ln(E(t))$ as a function of $t$ and we overlay a line defined by $ -1.2t-1.85$, this confirms the exponential decay rate of the energy, precisely we have 
$$E(t)\leq C e^{-1.2t}, \ \text{with} \ C=e^{-1.85}.$$
\begin{figure}[!ht]
\includegraphics[width=6.5cm]{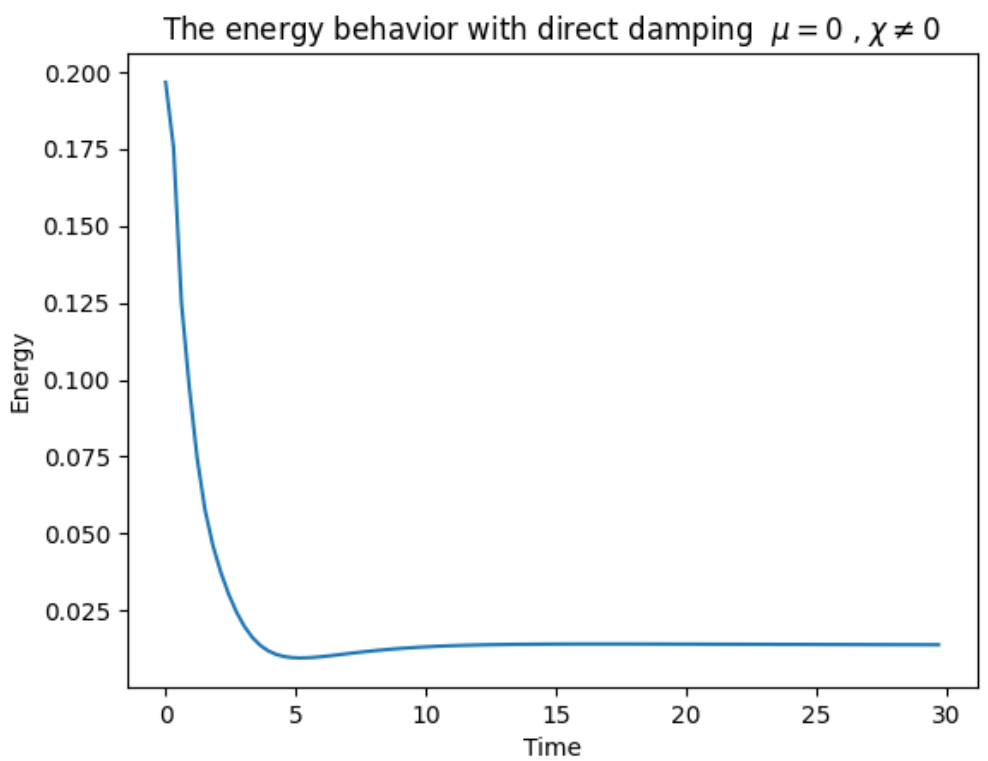}
\includegraphics[width=6.5cm]{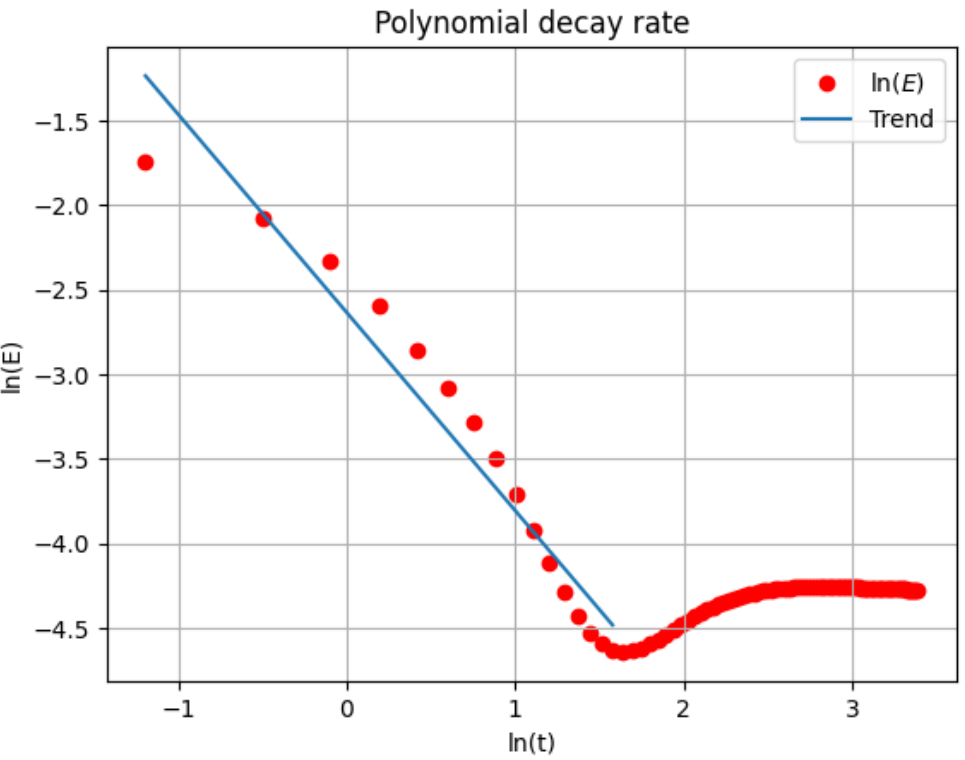}
\centering
\caption{The asymptotic behavior of the energy, $\mu=0$ and $\chi\neq 0$ .}
\label{fig:ln}
\end{figure}
 When the stability number is not equal to zero, the solution has a polynomial decay rate, as shown in the final time profile of the energy on the left-hand side of Figure \ref{fig:ln}. From the plot on the right side of Figure \ref{fig:ln}, we can write $\ln(E)\leq a \ln(t) +b$.  Based on our numerical analysis, we have derived an estimated value of the degree of the polynomial governing the decay rate associated with the discrete energy. Additionally, we observe that the discrete energy is approaching a non-zero constant at the final time $T=30$. This is in good agreement with the theoretical results in \cite{2,1}.
\subsubsection{Linear damping $\mu \psi_t$.}
\begin{figure}[!ht]
\includegraphics[width=6.5cm]{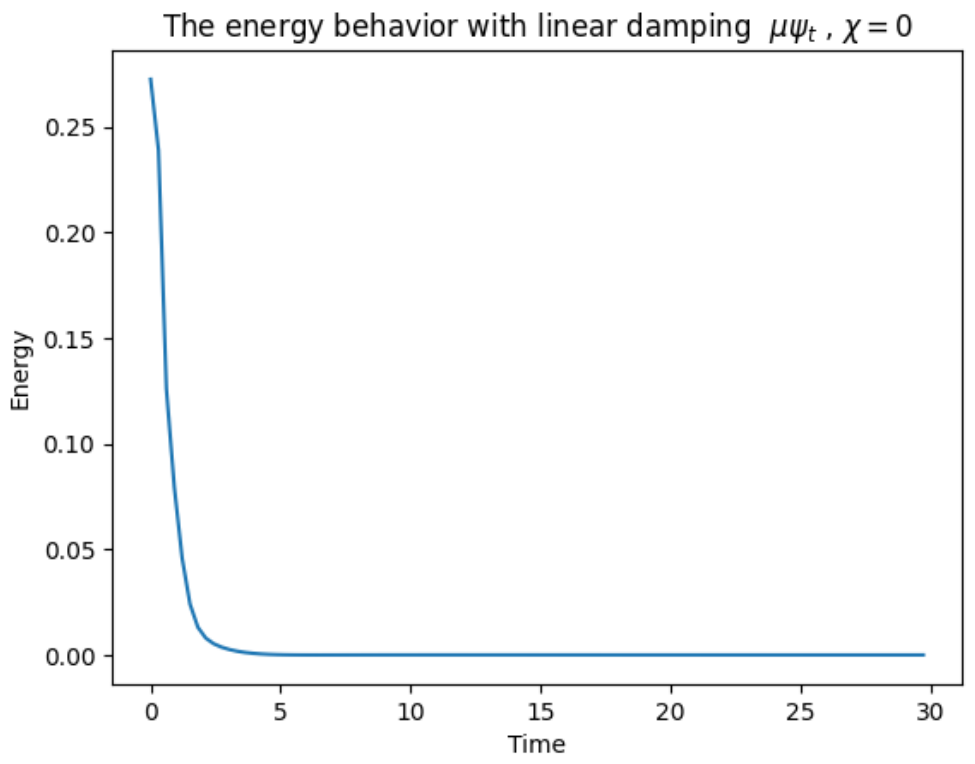}
\includegraphics[width=6.5cm]{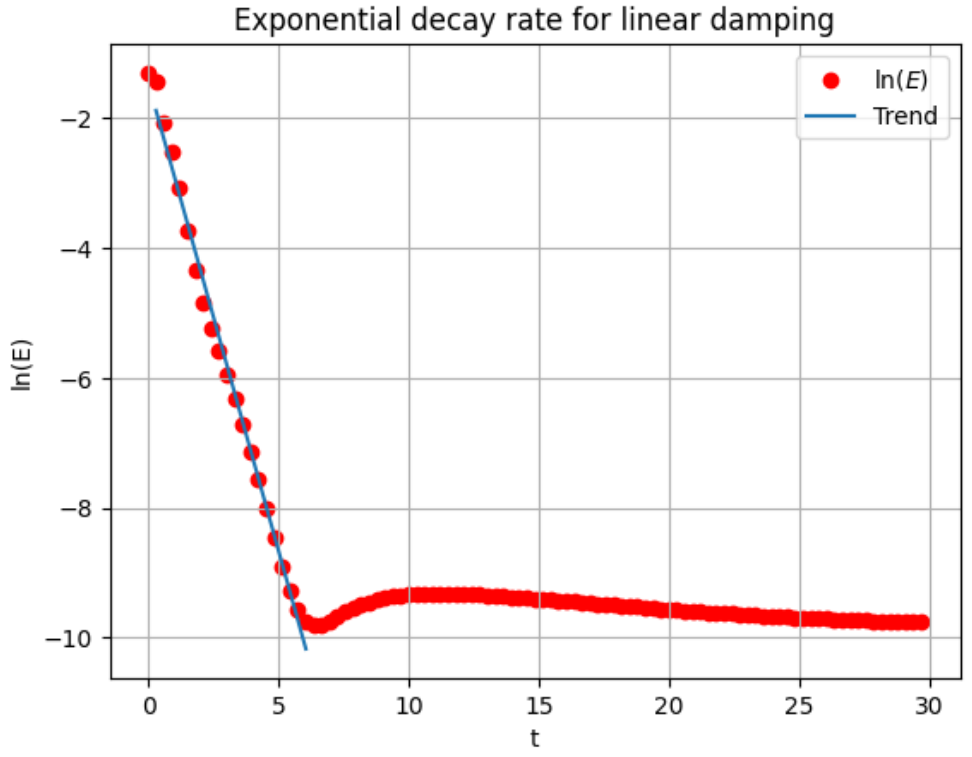}
\centering
\caption{The  exponential behavior of the energy when $\chi= 0$ and for a linear damping $\mu \psi_t$.}
\label{fig:lnll}
\end{figure}
Taking $\mu=1$ and $\chi=0$, Figure \ref{fig:lnll} shows decreasing energy over an extended period, with the right side of the graph illustrating the exponential decay rate of the energy. This decay rate is represented by the natural logarithm of energy, $\ln(E(t))$, plotted as a function of time $t$. The graph initially displays a decreasing linear trend from $t=0$ to $t=  5$, which is followed by a nearly constant value of almost equal to $-8.5$. As it is clearly shown in Figure \ref{fig:lnll}, the energy quickly attains a stable state equal to zero.
        \begin{figure}[!ht] 
        \centering
        \includegraphics[width=6.5cm]{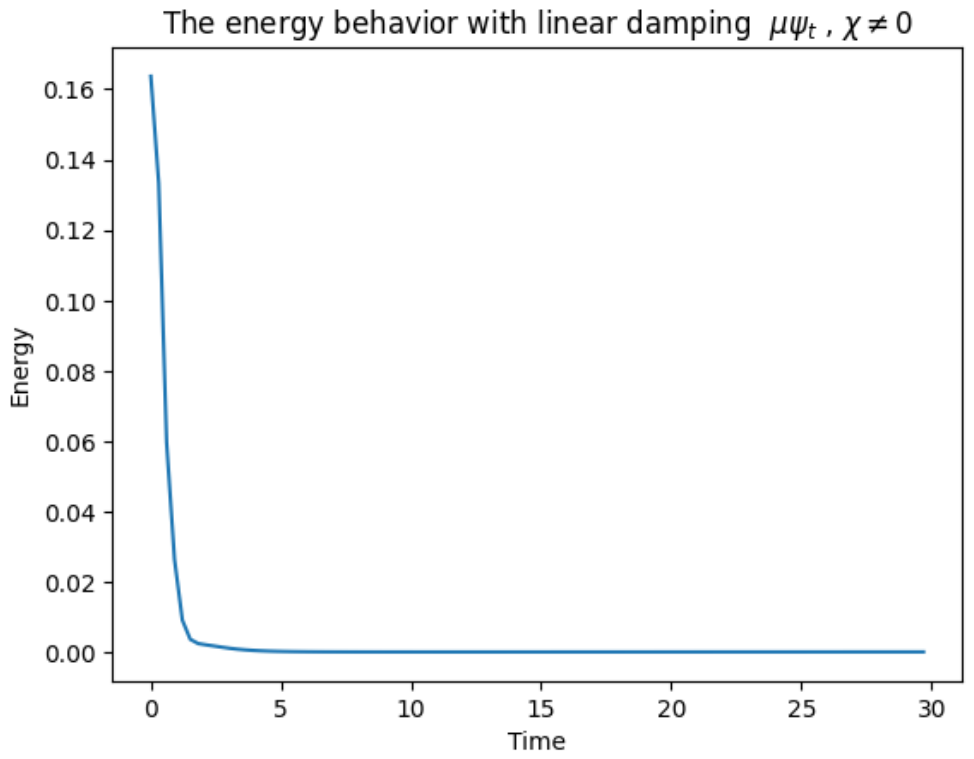}
        \includegraphics[width=6.5cm]{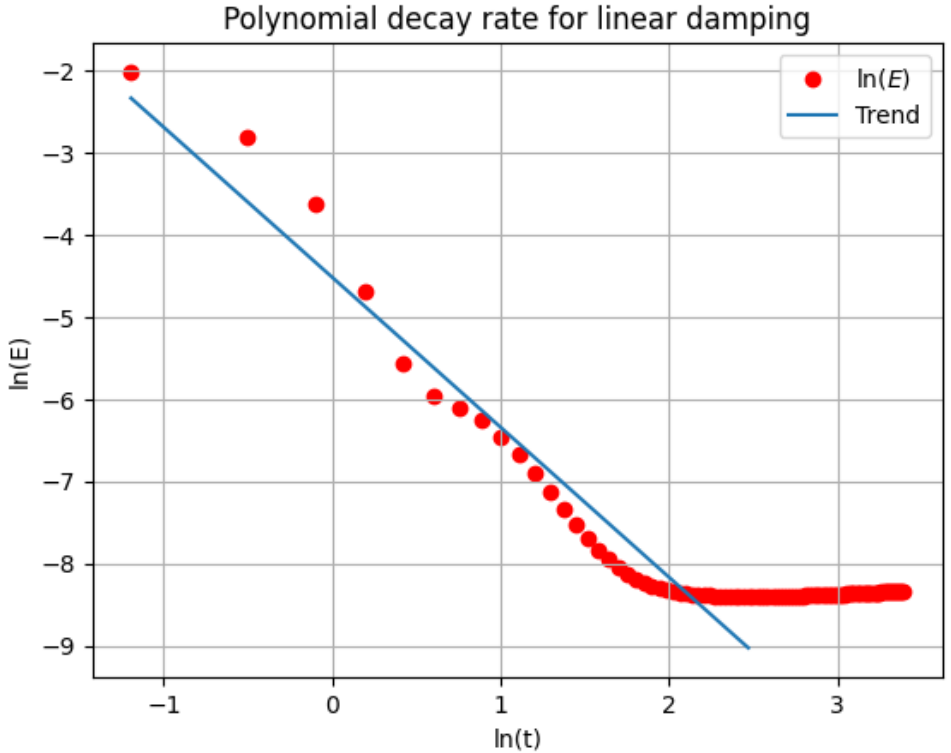}
        \caption{The polynomial decay of the energy when $\chi \neq 0$ and with the linear damping $\psi_t$.}
        \label{fig:lnl34}
        \end{figure}
    Figure \ref{fig:lnl34} shows that, when $\mu=1$ and $ \chi \neq 0$, the energy has a polynomial decay rate. This is shown on the right-hand side plot which exhibits the graph  of $\ln(E(t))$ in terms of $\ln(t)$.
\subsubsection{Non linear damping.}
   It is crucial to note that the shift from linear  to nonlinear damping requires the change of  the expression of the residuals associated with the nonlinearly damped Timoshenko system, {\textit i.e.}, replacing the linear damping $\mu \psi_t$ on the rotation equation of the system \eqref{1} by a nonlinear damping term of the form $\mu (\psi_t)^2$   and  $ \frac{\mu}{\psi_t} \cdot  \exp\left(\frac{-1}{\psi_t^2}\right)$, respectively.  The existence and uniqueness of \eqref{1} remain true since this has been proven in  \cite{1} in a more general context including  the above mentioned  cases. 
   
Therefore, we denote by $\Tilde{\mathcal{R}}_{PDE_2}$, $\Tilde{\Tilde{\mathcal{R}}}_{PDE_2}$ the new residuals corresponding to the second equation \eqref{1}$_2$, in which we impose   damping terms of the form $\mu (\psi_t)^2$ and $\frac{\mu}{\psi_t} \cdot  \exp\left(\frac{-1}{\psi_t^2}\right)$, respectively,  are given by
   \begin{equation} \label{Rnew}
\begin{array}{ll}
\tilde{\mathcal{R}}_{PDE_2}([\Tilde{\psi},\Tilde{\varphi},\Tilde{\theta}])(x,t)= \rho _{2}\Tilde{\psi} _{tt}-b\Tilde{\psi} _{xx}+k(\Tilde{\varphi} _{x}+\Tilde{\psi} )+\delta \Tilde{\theta
}_{x}+\mu(\tilde{\psi}_t)^2,
 \end{array}
  \end{equation}
  and 
   \begin{equation} \label{RRnew}
\begin{array}{ll}
\Tilde{\Tilde{\mathcal{R}}}_{PDE_2}([\Tilde{\psi},\Tilde{\varphi},\Tilde{\theta}])(x,t)= \rho _{2}\Tilde{\psi} _{tt}-b\Tilde{\psi} _{xx}+k(\Tilde{\varphi} _{x}+\Tilde{\psi} )+\delta \Tilde{\theta
}_{x}+\frac{\mu}{\tilde{\psi}_t} \cdot  \exp\left(\frac{-1}{\tilde{\psi}_t^2}\right).
 \end{array}
\end{equation}
The new training loss functions, as in  \eqref{LOSS}, will be calculated based on the new modified residuals \eqref{Rnew} and \eqref{RRnew}. 
\begin{figure}[!ht]
\includegraphics[width=6.5cm]{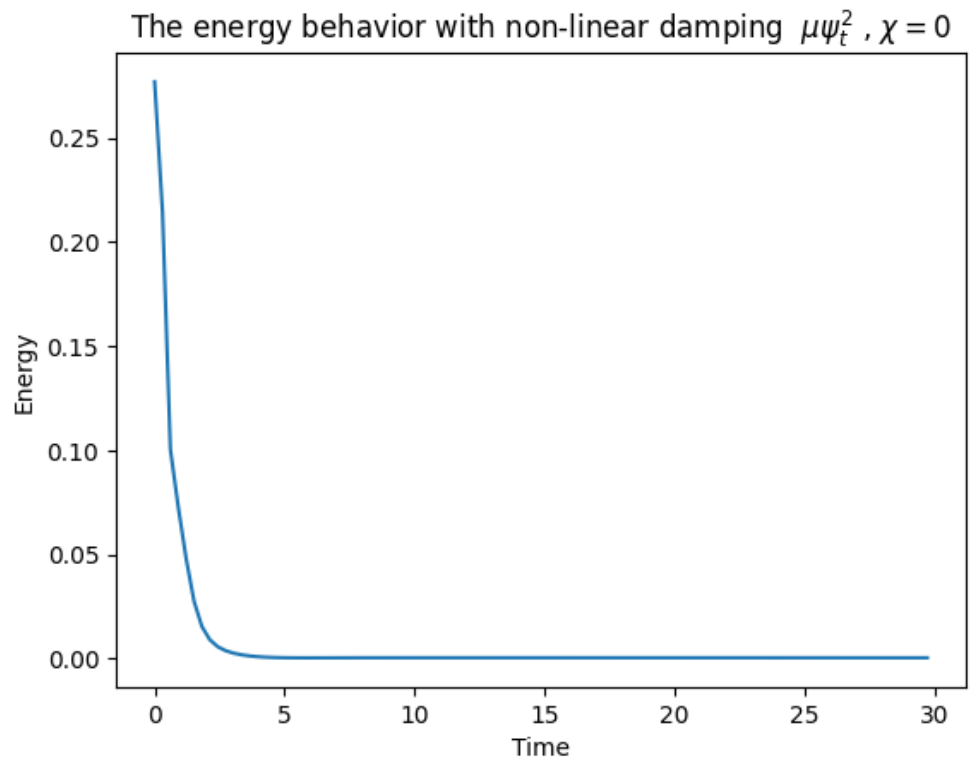}
\includegraphics[width=6.5cm]{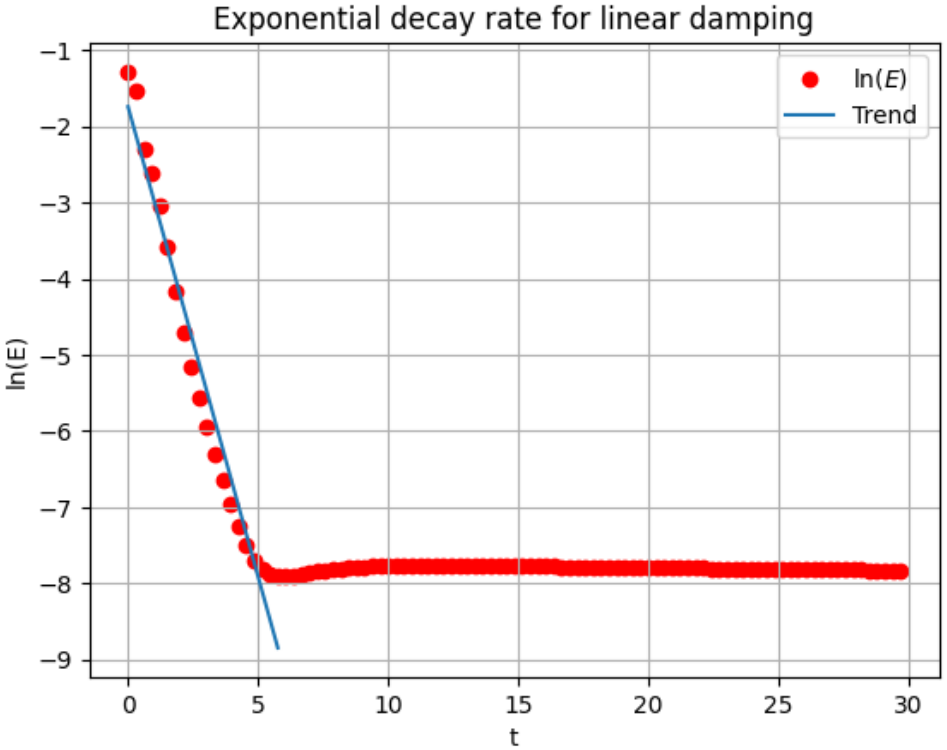}
\label{fig:lnl12}
\centering 
\caption{Graph for final  energy exponential behavior when $\chi= 0$ for a non-linearly damped system $\mu \psi_t^2$.}
\end{figure}

\begin{figure}[!ht]
\includegraphics[width=6.5cm]{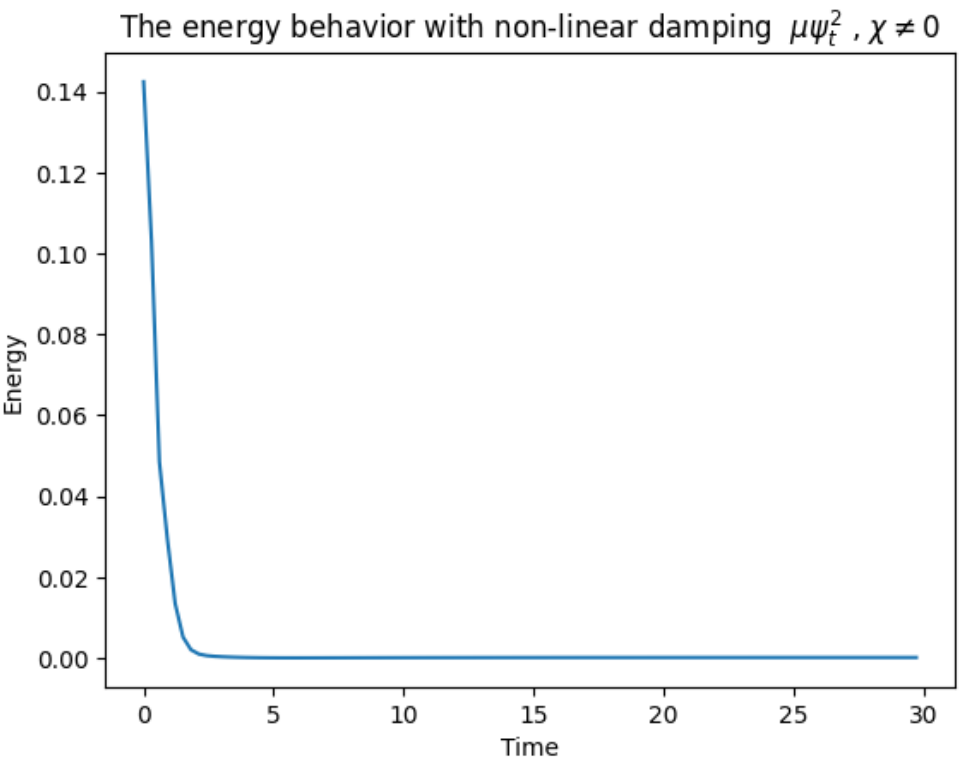}
 \includegraphics[width=6.5cm]{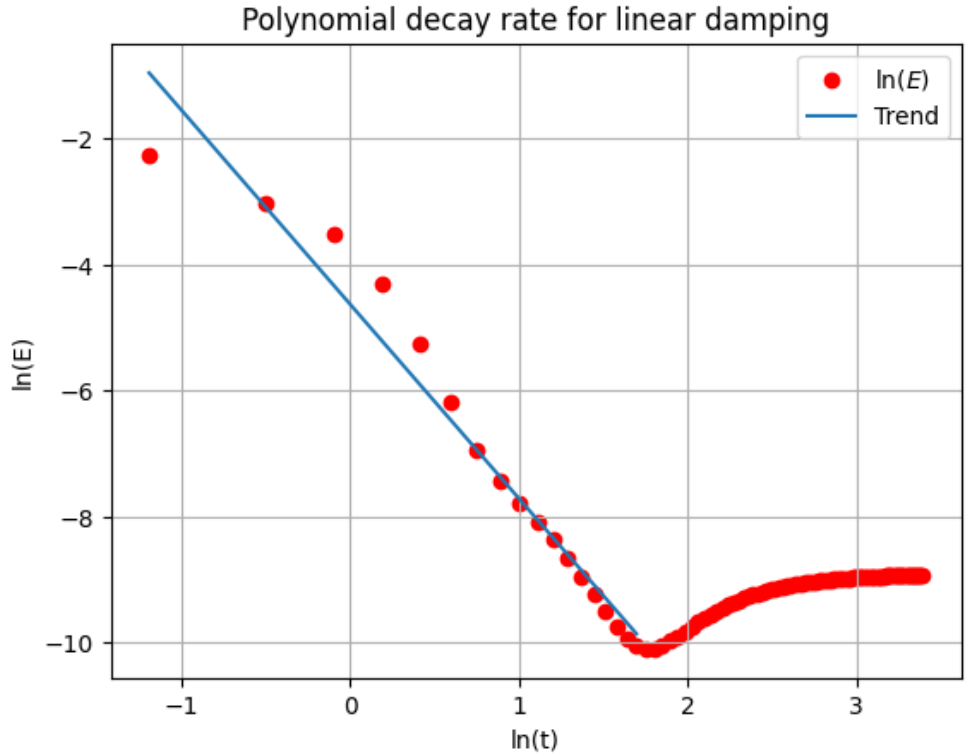}
\centering
\caption{Graph for final  energy polynomial behavior when $\chi \neq 0$ for a non-linearly damped system $\mu \psi_t^2$.}
\end{figure}
\subsubsection{Optimal logarithmic decay rate 
$\frac{\mu}{\psi_t} \cdot  \exp\left(\frac{-1}{\psi_t^2}\right)$.}
\begin{figure}[!ht]
\includegraphics[width=6.5cm]{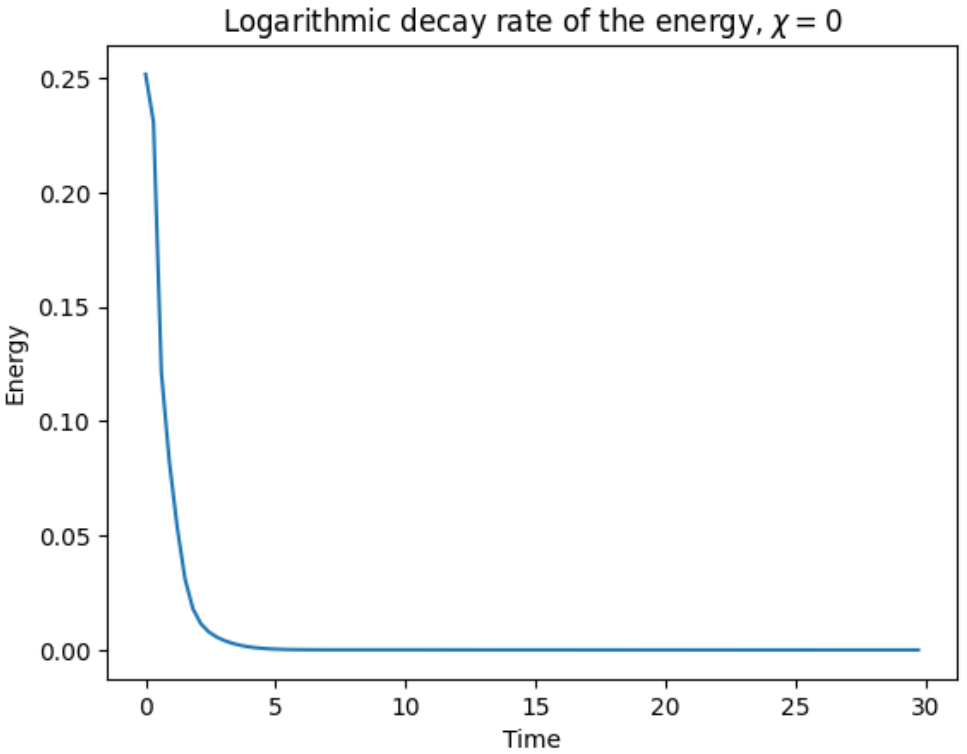}
\includegraphics[width=6.5cm]{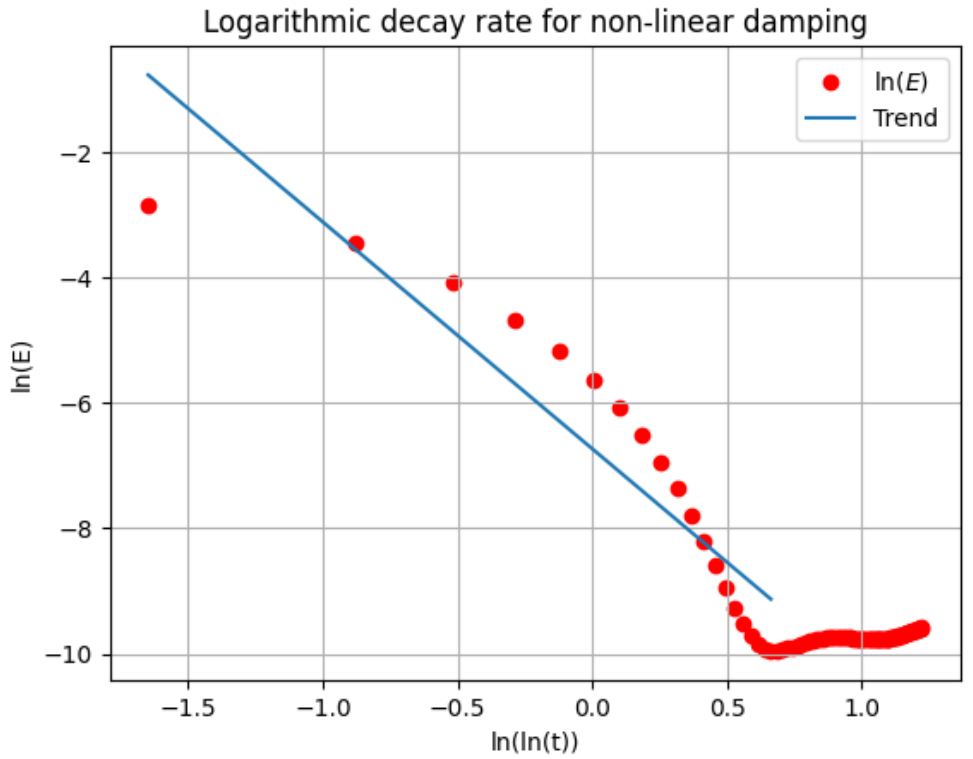}
\centering
\caption{The energy behavior when $\chi=0$.}
\label{fig:mlog2}
\end{figure}

We emphasize a specific example for the damping term $\frac{\mu}{\psi_t} \cdot  \exp\left(\frac{-1}{\psi_t^2}\right)$ regarding the energy associated with system \eqref{1}. This damping term  satisfies, for $t$ large enough, the following decay rate
\begin{equation}
c^{\prime }\left( \ln(t)\right)^{-1} \leq E(t)\leq c \left(
\ln(t)\right)^{-1}.
\end{equation}
The above bounds yield a more accurate energy decay rate.\\
To numerically test the theoretical result, we use the following parameters: $T=30$, $\chi=0$, $N_x=N_t=3000$. From this,  we see an asymptotic behavior of the energy as shown in Figure \ref{fig:mlog2}. The initial value of the energy is $\approx 0.28$. Then, the plot shows the evolution in time of the energy until getting to the stable value zero upon  $t \ge 2.75$.
Conversely, when $\chi\neq 0$, the initial value of the energy is $\approx 0.25$, and starting from $t=5$, the energy converges slowly to zero as shown in Figure \ref{fig:mlog2}. This is due to the weak nonlinear damping term added to the second equation. The right hand side of Figure \ref{fig:mlog2} shows that, by plotting $\ln(E(t))$ as a function of $\ln\left(\ln(t)\right)$,  the slow decay rate of the energy is logarithmic. 
\begin{figure}[!ht]
\includegraphics[width=6.5cm]{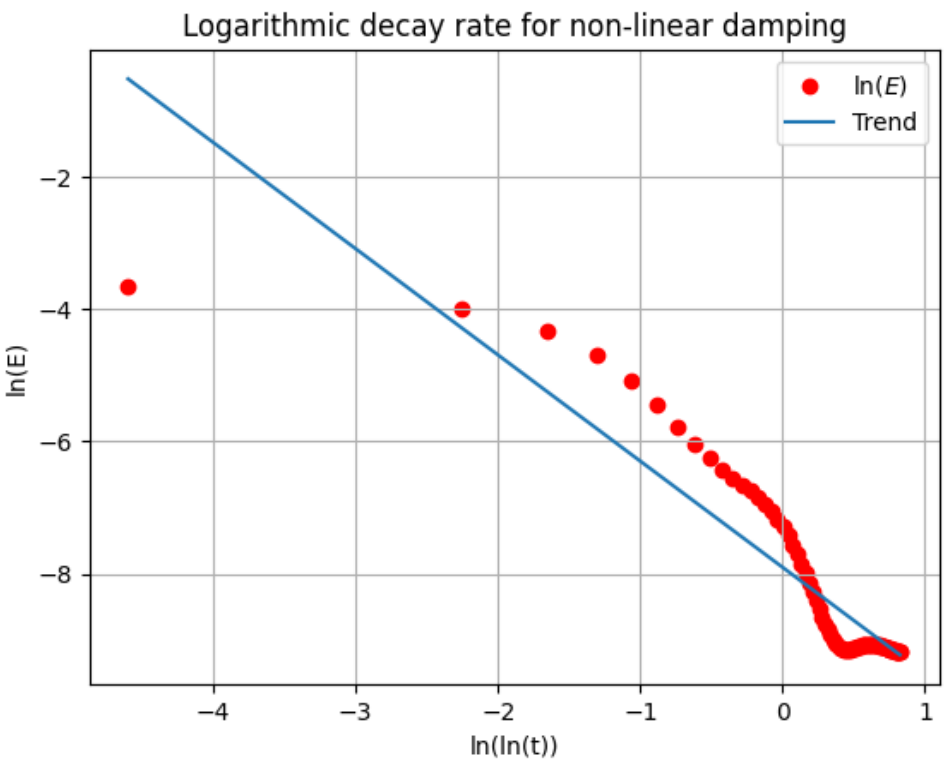}
\centering
\label{fig:og}
\caption{Logarithmic decay rate of the energy.}
\end{figure}
\section{Conclusion and future work}\label{sec4}
In this study, we considered  the Timoshenko system with thermoelasticity which represents a more comprehensive variant of the classical Timoshenko beam model, as it accounts for temperature and heat flux in addition to the beam's displacement and rotational angles during vibrations. Using the framework of Physics-Informed Neural Networks (PINNs), we could  approximate the Timoshenko solution to  \eqref{1} in this complex setting.  This achievement carries significant implications for addressing complex problems within the realms of applied mathematics and engineering using PINNs. We also determinated  the decay rate of discrete energy  \eqref{Eee} while considering a linear and a non-linear type  of dissipative mechanisms. As we are already aware, the study of Timoshenko system, when coupled with thermoelasticity, involves a stability number, denoted by $\chi$. Our numerical results unequivocally demonstrate that the exponential decay is achieved when $\chi=0$ without the addition of any dissipative terms. Therefore, we have examined the behavior of the energy when the system is linearly damped by introducing a damping term in the second equation of the system, specifically $\mu \psi_t$ and when the system is non linearly damped by imposing the term $\mu \psi_t^2$ and  $\frac{\mu}{\psi_t} \cdot  \exp\left(\frac{-1}{\psi_t^2}\right)$ to the rotation angle equation. In doing so, we have derived explicit decay rates, both exponential and polynomial as well as a logarithmic optimal decay rate for the discrete energy. One of the notable findings from our research is the ability to discern distinct energy decay rates when the system is subjected to various types of damping terms. This insight enhances our understanding of the system's dynamic behavior and contributes to the broader field of applied mathematics and engineering.\\
As we continue our exploration in this direction. We recall that while the Physics-Informed Neural Network (PINN) has found extensive applications in the approximation of partial differential equations (PDEs), the extent of theoretical exploration concerning its convergence and errors remains relatively restricted. Considering various facets of the neural network error is crucial, including factors like optimization error, approximation error, and estimation error \cite{0,Yeon}, primarily focused on elliptic and parabolic PDEs. Notably, there is a prominent dearth of comprehensive theoretical investigations into the convergence of PINN for coupled hyperbolic PDEs.











\section*{\bf\Large Acknowledgment}
The first author has been supported by the DAAD project n1-700136-01-EF, \textit{ Non Negative Structured Regression with Application in Communication and Data Science.}
Further, the first author would like to announce the non-financial support from The African Institute for Mathematical Sciences (AIMS), South Africa in the framework of the scientific collaboration with Technical University of Berlin.
\medskip

\hrule  

\end{document}